\documentclass[12pt,a4paper,twoside]{article}

\newcommand{\pageformat}[6]{\setlength{\hoffset}{-1in}
                  \setlength{\voffset}{-1in}
                  \addtolength{\hoffset}{#5}
                            \addtolength{\voffset}{#6}
                            \setlength{\oddsidemargin}{#1}
                            \setlength{\evensidemargin}{#2}
                            \setlength{\textwidth}{\paperwidth}
                  \addtolength{\textwidth}{-\oddsidemargin}
                  \addtolength{\textwidth}{-\evensidemargin}
                  \addtolength{\textwidth}{-\marginparsep}
                  \addtolength{\textwidth}{-\marginparwidth}
                            \setlength{\topmargin}{#3}
                            \setlength{\textheight}{\paperheight}
                  \addtolength{\textheight}{-\topmargin}
                  \addtolength{\textheight}{-\headheight}
                  \addtolength{\textheight}{-\headsep}
                  \addtolength{\textheight}{-\footskip}
                  \addtolength{\textheight}{-#4}}
\pageformat{2cm}{3cm}{25mm}{25mm}{1pt}{0pt}

\usepackage{ifthen}
\newboolean{article}
    \setboolean{article}{true}
\newboolean{report}
\newboolean{book}
\newboolean{letter}
\newboolean{german}
\newboolean{italian}
\newboolean{nobaselinestretch}
\newboolean{nosectionappendix}
\newboolean{oldtoc}
\newboolean{nosectionequn}
\newboolean{notheorem}
\ifthenelse{\boolean{german}}{
    \usepackage{german}}{}

\usepackage[latin1]{inputenc}
\usepackage{amsmath}
\usepackage{amssymb}
\usepackage[mathscr]{eucal}

\ifthenelse{\boolean{notheorem}}{}{
    \usepackage{theorem}}

\ifthenelse{\boolean{nobaselinestretch}}{}{
    \renewcommand{\baselinestretch}{1.25}}

\newenvironment{env}[2]{\begin{#1}#2\end{#1}}{}
    \newcommand{\beq}[1]{\begin{env}{equation}{#1}}
    \newcommand{\beqn}[1]{\begin{env}{equation*}{#1}}
    \newcommand{\bal}[1]{\begin{env}{align}{#1}}
    \newcommand{\baln}[1]{\begin{env}{align*}{#1}}
    \newcommand{\bga}[1]{\begin{env}{gather}{#1}}
    \newcommand{\bgan}[1]{\begin{env}{gather*}{#1}}
    \newcommand{\bflal}[1]{\begin{env}{flalign}{#1}}
    \newcommand{\bflaln}[1]{\begin{env}{flalign*}{#1}}
    \newcommand{\bmu}[1]{\begin{env}{multline}{#1}}
    \newcommand{\bmun}[1]{\begin{env}{multline*}{#1}}
    \newcommand{\bsp}[1]{\begin{env}{split}{#1}}

    \newcommand{\eeq}{\end{env}}
    \newcommand{\eeqn}{\end{env}}
    \newcommand{\eal}{\end{env}}
    \newcommand{\ealn}{\end{env}}
    \newcommand{\ega}{\end{env}}
    \newcommand{\egan}{\end{env}}
    \newcommand{\eflal}{\end{env}}
    \newcommand{\eflaln}{\end{env}}
    \newcommand{\emu}{\end{env}}
    \newcommand{\emun}{\end{env}}
    \newcommand{\esp}{\end{env}}

\newcommand{\lf}{\vspace{2ex}}

\renewcommand{\bf}[1]{\textbf{#1}}
\renewcommand{\it}[1]{\textit{#1}}

\renewcommand{\tt}[1]{\texttt{#1}}
\newcommand{\hl}[1]{\bf{\it{#1}}}

\newcommand{\cmc}[1]{\mathcal{#1}}
\newcommand{\eus}[1]{\mathscr{#1}}

\newcommand{\bb}[1]{\mathbb{#1}}

\newcommand{\nbd}[1]{$#1$\nobreakdash--}
\newcommand{\ol}[1]{\overline{#1}}

\newcommand{\abs}[1]{\left\lvert#1\right\rvert}
\newcommand{\norm}[1]{\left\lVert#1\right\rVert}
\newcommand{\babs}[1]{\bigl\lvert#1\bigr\rvert}

\newcommand{\sabs}[1]{\abs{\smash{#1}}}

\newcommand{\bfam}[1]{\bigl(#1\bigr)}

\newcommand{\AB}[1]{\langle#1\rangle}

\newcommand{\BAB}[1]{\Bigl\langle#1\Bigr\rangle}
\newcommand{\CB}[1]{\{#1\}}
\newcommand{\bCB}[1]{\bigl\{#1\bigr\}}
\newcommand{\BCB}[1]{\Bigl\{#1\Bigr\}}
\newcommand{\SB}[1]{[#1]}

\newcommand{\set}[2][]{
    \ifthenelse{\equal{#1}{}}{
        \CB{#2}}{
        \CB{#1~|~#2}}}
\newcommand{\bset}[2][]{
    \ifthenelse{\equal{#1}{}}{
        \bCB{#2}}{
        \bCB{#1~|~#2}}}
\newcommand{\Bset}[2][]{
    \ifthenelse{\equal{#1}{}}{
        \BCB{#2}}{
        \BCB{#1~\big|~#2}}}

\newcommand{\N}{\bb{N}}

\newcommand{\R}{\bb{R}}

\newcommand{\cA}{\cmc{A}}

\newcommand{\cL}{\cmc{L}}

\newcommand{\sL}{\eus{L}}

\newcommand{\I}{{I\!\!\!\;I}}

\ifthenelse{\boolean{nosectionequn}}{}{
    \numberwithin{equation}{section}
    }

\ifthenelse{\boolean{article}\or\boolean{letter}\or\boolean{nosectionequn}}{
    \setboolean{nosectionappendix}{true}}{}
\ifthenelse{\boolean{nosectionappendix}}{}{
    \renewcommand{\appendix}{
        \chapter*{\appendixname}
        \addcontentsline{toc}{chapter}{\appendixname}
        \renewcommand{\thesection}{\Alph{section}}
        \setcounter{section}{0}}}
   
\ifthenelse{\boolean{report}\or\boolean{book}}{
    }{}

\ifthenelse{\boolean{notheorem}}{}{
        \newcommand{\mnname}{Mathematical note.}
        \newcommand{\enname}{End of the note.}
        \newcommand{\definame}{Definition.}
        \newcommand{\propname}{Proposition.}
        \newcommand{\lemname}{Lemma.}
        \newcommand{\exname}{Example.}
        \newcommand{\exername}{Exercise.}
        \newcommand{\remname}{Remark.}
        \newcommand{\obname}{Observation.}
        \newcommand{\thmname}{Theorem.}
        \newcommand{\corname}{Corollary.}
        \newcommand{\proofname}{Proof.}
    \ifthenelse{\boolean{german}}{
        \renewcommand{\mnname}{Mathematische Notiz.}
        \renewcommand{\enname}{Ende der Notiz.}
        \renewcommand{\exname}{Beispiel.}
        \renewcommand{\exername}{Übung.}
        \renewcommand{\remname}{Bemerkung.}
        \renewcommand{\obname}{Beobachtung.}
        \renewcommand{\thmname}{Satz.}
        \renewcommand{\corname}{Korollar.}
        \renewcommand{\proofname}{Beweis.}}{}
    \ifthenelse{\boolean{italian}}{
        \renewcommand{\mnname}{Nota matematica.}
        \renewcommand{\enname}{Fina della nota.}
        \renewcommand{\definame}{Definizione.}
        \renewcommand{\propname}{Proposizione.}
        \renewcommand{\exname}{Esempio.}
        \renewcommand{\exername}{Esercizio.}
        \renewcommand{\remname}{Nota.}
        \renewcommand{\obname}{Osservazione.}
        \renewcommand{\thmname}{Teorema.}
        \renewcommand{\corname}{Corollario.}
        \renewcommand{\proofname}{Dimostrazione.}

       \renewcommand{\appendixname}{Appendice}

       }{}
    \theoremheaderfont{\normalfont\bfseries}
    \theoremstyle{change}
        \theorembodyfont{\rmfamily}
            \newtheorem{emp}{}[section]
                \newcommand{\bemp}[1][]{
                    \begin{emp}\hskip-\labelsep\bf{#1}\hskip\labelsep}
                \newcommand{\eemp}{\end{emp}}
\newtheorem{itemp}[emp]{}
                \newcommand{\bitemp}[1][]{
                    \begin{itemp}\hskip-\labelsep\bf{#1}\hskip\labelsep\normalfont\itshape}
                \newcommand{\eitemp}{\end{itemp}}
            \newtheorem{mn}[emp]{\mnname}
                \newcommand{\bnm}{\begin{mn}~\begin{quotation}\renewcommand{\baselinestretch}{1}\small\noindent\ignorespaces}
                \newcommand{\enm}{\end{quotation}\hfill\bf{\enname}\end{mn}}
            \newtheorem{ex}[emp]{\exname}
                \newcommand{\bex}{\begin{ex}}
                \newcommand{\eex}{\end{ex}}
            \newtheorem{exer}[emp]{\exername}
                \newcommand{\bexer}{\begin{exer}}
                \newcommand{\eexer}{\end{exer}}
            \newtheorem{defi}[emp]{\definame}
                \newcommand{\bdefi}{\begin{defi}}
                \newcommand{\edefi}{\end{defi}}
            \newtheorem{rem}[emp]{\remname}
                \newcommand{\brem}{\begin{rem}}
                \newcommand{\erem}{\end{rem}}
            \newtheorem{ob}[emp]{\obname}
                \newcommand{\bob}{\begin{ob}}
                \newcommand{\eob}{\end{ob}}
        \theorembodyfont{\normalfont\itshape}
            \newtheorem{thm}[emp]{\thmname}
                \newcommand{\bthm}{\begin{thm}}
                \newcommand{\ethm}{\end{thm}}
            \newtheorem{prop}[emp]{\propname}
                \newcommand{\bprop}{\begin{prop}}
                \newcommand{\eprop}{\end{prop}}
            \newtheorem{cor}[emp]{\corname}
                \newcommand{\bcor}{\begin{cor}}
                \newcommand{\ecor}{\end{cor}}
            \newtheorem{lem}[emp]{\lemname}
                \newcommand{\blem}{\begin{lem}}
                \newcommand{\elem}{\end{lem}}
\newenvironment{empn}[1]{\lf\noindent\bf{#1}\ignorespaces\hskip\labelsep}{\lf}
		\newcommand{\bempn}[1]{\begin{empn}{#1}}
		\newcommand{\eempn}{\end{empn}}
		\newcommand{\bitempn}[1]{\begin{empn}{#1}\normalfont\itshape}
		\newcommand{\eitempn}{\end{empn}}
                \newcommand{\bnmn}{\begin{empn}{\mnname}~\begin{quotation}\renewcommand{\baselinestretch}{1}\small\noindent\ignorespaces}
                \newcommand{\enmn}{\end{quotation}\hfill\bf{\enname}\end{empn}}
		\newcommand{\bexn}{\begin{empn}{\exname}}
		\newcommand{\eexn}{\end{empn}}
		\newcommand{\bexern}{\begin{empn}{\exername}}
		\newcommand{\eexern}{\end{empn}}
		\newcommand{\bdefin}{\begin{empn}{\definame}}
		\newcommand{\edefin}{\end{empn}}
		\newcommand{\bremn}{\begin{empn}{\remname}}
		\newcommand{\eremn}{\end{empn}}
		\newcommand{\bobn}{\begin{empn}{\obname}}
		\newcommand{\eobn}{\end{empn}}

\newcommand{\qedsymbol}{~\rule[-0.35mm]{2mm}{2mm}}
    \newcounter{proof}[emp]
    \newenvironment{Proof}[1]{
        \vspace{1ex}
        \renewcommand{\item}[1][\stepcounter{proof}(\roman{proof})]%
            {##1\hskip\labelsep}
        \noindent\textsc{#1\hskip\labelsep}}{
        \nolinebreak\qedsymbol}
    \newcommand{\proof}[1][\proofname]{
        \begin{Proof}{#1}\ignorespaces}
    \newcommand{\qed}{\end{Proof}}
    \newcommand{\noqed}{
        \renewcommand{\qedsymbol}{}
        \end{Proof}}}
    \ifthenelse{\boolean{italian}}{
        \renewcommand{\proofname}{Dimostrazione.}}{}

\usepackage[varg]{txfonts}

\usepackage[hypertex]{hyperref}

\begin{document}

\bibliographystyle{amsalpha}

\title{Extending the Set of Quadratic Exponential Vectors\thanks{LA and MS are supported by Italian MUR (PRIN 2007). MS is supported by research funds of the Dipartimento S.E.G.e S.\ of University of Molise.}}

\author{Luigi Accardi, Ameur Dharhi, and Michael Skeide}

\date{November 2008}

\maketitle

%\vspace{-3ex}
\begin{abstract}
\noindent
We extend the square of white noise algebra over the step functions on $\R$ to the test function space $L^2(\R^d)\cap L^\infty(\R^d)$, and we show that in the Fock representation the exponential vectors exist for all test functions bounded by $\frac{1}{2}$.
\end{abstract}

\section{Introduction}

%\noindent
Modulo minor variations in the choice of the test function space, the square of white noise (SWN) algebra has been introduced by Accardi, Lu and Volovich \cite{ALV99p} as follows. Let $\sL=L^2(\R^d)\cap L^\infty(\R^d)$ and $c>0$ a constant. Then the \hl{SWN algebra} $\cA$ over $\sL$ is the unital \nbd{*}algebra generated by symbols $B_f,N_f$ $(f\in\sL)$ and the commutation relations
\baln{
\SB{B_f,B_g^*}
&
~=~
2c\AB{f,g}+4N_{\ol{f}g},
&
\SB{N_f,B_g^*}
&
~=~
2B_{fg}^*,
}\ealn
$(f,g\in\sL)$ and all other commutators $0$. Note that by the first relation, $N_f^*=N_{\ol{f}}$.

A \hl{Fock representation} of $\cA$ is a representation ($*$, of course) $\pi$ of $\cA$ on a pre-Hilbert space $H$ with a unit vector $\Phi\in H$, fulfilling $\cA\Phi=H$ and $\pi(B_f)\Phi=\pi(N_f)\Phi=0$ for all $f\in\sL$. From the commutation relations it follows that a Fock representation is unique up to unitary equivalence. Existence of a Fock representation has been established by different proofs in \cite{ALV99p,AcSk00a,Sni00,AFS02} for $d=1$. They extend easily to general $d\in\N$. Henceforth, we speak about \bf{the} Fock representation. The Fock representation would be faithful, if we require also that the $N_f$ depend linearly on $f$. By abuse of notation, we identify $\cA$ with its image $\pi(\cA)$ omitting, henceforth, $\pi$.

The \hl{exponential vector} $\psi(f)$ to an element $f\in\sL$ is defined as
\beqn{
\psi(f)
~:=~
\sum_{m=0}^\infty\frac{{B_f^*}^m\Phi}{m!}
}\eeqn
whenever the series exists. In Accardi and Skeide \cite{AcSk00} is has been shown for $d=1$ that $\psi(\sigma\I_{\SB{0,t}})$ exists for $\abs{\sigma}<\frac{1}{2}$ and that $\AB{\psi(\sigma\I_{\SB{0,t}}),\psi(\rho\I_{\SB{0,t}})}=e^{-\frac{ct}{2}\ln(1-4\ol{\sigma}\rho)}$. As noted in \cite{AcSk00}, this extends to arbitrary step functions $f,g$ on $\R$ with $\norm{f}_\infty<\frac{1}{2}$, with inner product
\beqn{\tag{$*$}\label{*}
\AB{\psi(f),\psi(g)}
~=~
e^{-\frac{c}{2}\int\ln(1-4\ol{f(t)}g(t))\,dt}.\footnote{The \it{correlation kernel} on the right-hand side coincides, modulo scaling, with the correlation kernel in Boukas' representation \cite{Bou91} of Feinsilver's \it{finite difference algebra} \cite{Fei87}. In \cite{AcSk00}, this observation gave rise to the discovery of an intimate relation between the SWN algebra and the finite difference algebra.}
}\eeqn
Our scope is to extend the set of exponential vectors and the formula in \eqref{*} for their inner product to test functions $f\in\cL$ with $\norm{f}_\infty<\frac{1}{2}$.

In the ``29th Quantum Probabililty Conference'' in October 2008 in Hammamet, Tunisia, Dhahri explained that the extension can be done for exponential vectors to all elements $f$ in $\sL$ with $\norm{f}_\infty<\frac{1}{2}$. This a part of the work Accardi and Dhahri \cite{AcDh08p} (in preparation) on the \it{second quantization functor} for the square of white noise. Here we give a simple proof of this partial result.

\section{The result}

\bthm
The exponential vector $\psi(f)$ exists for every $f\in\cL$ with $\norm{f}_\infty<\frac{1}{2}$ and the inner product of two such exponential vectors is given by \eqref{*}.
\ethm

\proof
\item
We show that the right-hand side of \eqref{*} exists. Indeed, by Taylor expansion we have $\abs{\ln(1+x)}\le M_\delta\abs{x}$ for $\abs{x}\le1-\delta$ for every $\delta\in(0,1)$, where $M_\delta$ may depend on $\delta$ but not on $x$. Choose $\delta=1-4\norm{f}_\infty\norm{g}_\infty\in(0,1)$. Then
\beqn{
\babs{\ln(1-4\ol{f(t)}g(t))}
~\le~
M_\delta\babs{4\ol{f(t)}g(t)}.
}\eeqn
Since $\babs{\ol{f(t)}g(t)}$ is integrable, so is $\ln(1-4\ol{f(t)}g(t))$.

\item
The function $x\mapsto\ln x$ is increasing on the whole half line $(0,\infty)$. It follows that also the function $x\mapsto-\ln(1-x)$ is increasing on $(-1,1)$. We conclude that $\frac{1}{2}>\abs{f}\ge\abs{g}$ implies $-\ln(1-4\abs{f(t)}^2)\ge-\ln(1-4\abs{g(t)}^2)$. Choose for $f$ an \nbd{L^2}approximating sequence of step functions $\bfam{f_n}_{n\in\N}$ in such a way that $\abs{f}\ge\abs{f_n}$ for all $n\in\N$. By the \it{dominated convergence theorem}, $\lim_{n\to\infty}e^{-\frac{c}{2}\int\ln(1-4\abs{f_n(t)}^2)\,dt}=e^{-\frac{c}{2}\int\ln(1-4\abs{f(t)}^2)\,dt}$.

\item
In precisely the same way as in \cite{AcSk00}, one shows that \eqref{*} is true for all step functions strictly bounded by $\frac{1}{2}$. It follows that $\lim_{n\to\infty}\norm{\psi(f_n)}^2=e^{-\frac{c}{2}\int\ln(1-4\abs{f(t)}^2)\,dt}$.

\item
Since $\AB{{B_f^*}^m\Phi,{B_f^*}^m\Phi}$ is a polinomial (of degree $m$) in $\AB{f,f}$, it depends continuously in \nbd{L^2}norm on $f$. So, for every $M\in\N$ there is an $n\in\N$ such that
\bmun{
\BAB{\sum_{m=0}^M\frac{{B_f^*}^m\Phi}{m!},\sum_{m=0}^M\frac{{B_f^*}^m\Phi}{m!}}
~\le~
\BAB{\sum_{m=0}^M\frac{{B_{f_n}^*}^m\Phi}{m!},\sum_{m=0}^M\frac{{B_{f_n}^*}^m\Phi}{m!}}+1
\\
~\le~
\BAB{\sum_{m=0}^\infty\frac{{B_{f_n}^*}^m\Phi}{m!},\sum_{m=0}^\infty\frac{{B_{f_n}^*}^m\Phi}{m!}}+1
~=~
\norm{\psi(f_n)}^2+1
~\le~
e^{-\frac{c}{2}\int\ln(1-4\abs{f(t)}^2)\,dt}+1.
}\emun
By the theorem on exchange of limits under domination, it follows that
\bmun{
\lim_{M\to\infty}\BAB{\sum_{m=0}^M\frac{{B_f^*}^m\Phi}{m!},\sum_{m=0}^M\frac{{B_f^*}^m\Phi}{m!}}
~=~
\lim_{M\to\infty}\lim_{n\to\infty}\BAB{\sum_{m=0}^M\frac{{B_{f_n}^*}^m\Phi}{m!},\sum_{m=0}^M\frac{{B_{f_n}^*}^m\Phi}{m!}}
\\
~=~
\lim_{n\to\infty}\lim_{M\to\infty}\BAB{\sum_{m=0}^M\frac{{B_{f_n}^*}^m\Phi}{m!},\sum_{m=0}^M\frac{{B_{f_n}^*}^m\Phi}{m!}}
~=~
\lim_{n\to\infty}\norm{\psi(f_n)}^2
~=~
e^{-\frac{c}{2}\int\ln(1-4\abs{f(t)}^2)\,dt}.
}\emun
From this we conlcude that $\psi(f)$ exists and that $\norm{\psi(f)}^2=e^{-\frac{c}{2}\int\ln(1-4\abs{f(t)}^2)\,dt}$.

\item
Doing the same sort of computation for the difference $\psi(f)-\psi(f_n)$, it follows that $\lim_{n\to\infty}\psi(f_n)=\psi(f)$. Approximating also $g$ by a sequence of step functions $g_n$ with $\abs{g}\ge\abs{g_n}$, we find $\lim_{n\to\infty}\AB{\psi(f_n),\psi(g_n)}=\AB{\psi(f),\psi(g)}$ (continuity of the inner product), and
\beqn{
\lim_{n\to\infty}e^{-\frac{c}{2}\int\ln(1-4\ol{f_n(t)}g_n(t))\,dt}
~=~
e^{-\frac{c}{2}\int\ln(1-4\ol{f(t)}g(t))\,dt}
}\eeqn
(once more, by dominated convergence for $\sabs{\ol{f_n}g_n}\le\sabs{\ol{f}g}$ on the other side. This shows \eqref{*} for all $f,g$ as specified.\qed

\setlength{\baselineskip}{2.5ex}

%\bibliography{mybib}

\begin{thebibliography}{ALV99}

\bibitem[AD08]{AcDh08p}
L.~Accardi and A.~Dhahri, \emph{{Second quantization on the square of white
  noise algebra}}, Pre\-print, Rome (in preparation), 2008.

\bibitem[AFS02]{AFS02}
L.~Accardi, U.~Franz, and M.~Skeide, \emph{{Renormalized squares of white noise
  and other non-gau{\ss}ian noises as L\'evy processes on real Lie algebras}},
  Commun.\ Math.\ Phys. \textbf{228} (2002), 123--150, (Rome,
  Volterra-Pre\-print 2000/0423).

\bibitem[ALV99]{ALV99p}
L.~Accardi, Y.G. Lu, and I.V. Volovich, \emph{{White noise approach to
  classical and quantum stochastic calculi}}, Pre\-print, Rome, 1999, To appear
  in the lecture notes of the Volterra International School of the same title,
  held in Trento.

\bibitem[AS00a]{AcSk00a}
L.~Accardi and M.~Skeide, \emph{{Hilbert module realization of the square of
  white noise and the finite difference algebra}}, Math.\ Notes \textbf{86}
  (2000), 803--818, (Rome, Volterra-Pre\-print 1999/0384).

\bibitem[AS00b]{AcSk00}
\bysame, \emph{{On the relation of the square of white noise and the finite
  difference algebra}}, Infin.\ Dimens.\ Anal.\ Quantum Probab.\ Relat.\ Top.
  \textbf{3} (2000), 185--189, (Rome, Volterra-Pre\-print 1999/0386).

\bibitem[Bou91]{Bou91}
A.~Boukas, \emph{{An example of quantum exponential process}}, Mh.\ Math.
  \textbf{112} (1991), 209--215.

\bibitem[Fei87]{Fei87}
P.J. Feinsilver, \emph{{Discrete analogues of the Heisenberg-Weyl algebra}},
  Mh.\ Math. \textbf{104} (1987), 89--108.

\bibitem[Sni00]{Sni00}
P.~Sniady, \emph{{Quadratic bosonic and free white noise}}, Commun.\ Math.\
  Phys. \textbf{211} (2000), 615--628.

\end{thebibliography}
\newcommand{\Swap}[2]{#2#1}\newcommand{\Sort}[1]{}
\providecommand{\bysame}{\leavevmode\hbox to3em{\hrulefill}\thinspace}
\providecommand{\MR}{\relax\ifhmode\unskip\space\fi MR }
\providecommand{\MRhref}[2]{%
  \href{http://www.ams.org/mathscinet-getitem?mr=#1}{#2}
}
\providecommand{\href}[2]{#2}

\noindent
Luigi Accardi: \it{Centro Vito Volterra, Facoltà di Economia, Università degli Studi di Roma ``Tor Vergata'', Via Columbia 2, 00133 Roma, Italy},
E-mail: \tt{accardi@volterra.uniroma2.it},\\
Homepage: \tt{http://www.volterra.uniroma2.it}

\lf\noindent
Ameur Dhahri: \it{Ceremade, Université Paris Dauphine, Place de Lattre de Tassigny, 75775 Paris Cédex 16, France},
E-mail: \tt{dhahri@ceremade.dauphine.fr}

\lf\noindent
Michael Skeide: \it{Dipartimento S.E.G.e S., Università degli Studi del Molise, Via de Sanctis, 86100 Campobasso, Italy},
E-mail: \tt{skeide@math.tu-cottbus.de},\\
Homepage: \tt{http://www.math.tu-cottbus.de/INSTITUT/lswas/\_skeide.html}

\end{document}